\documentclass[]{article}
\usepackage{amssymb}

\newcommand{\dfrac}{\displaystyle \frac}
\newcommand{\dsum}{\displaystyle \sum}
\newcommand{\dint}{\displaystyle \int}

\begin{document}

\title{On DSS Method for a Fast Identification of the Static and Dynamic Responses
of Servovalves\thanks{%
This research was supported by the Russian Found for Basic research, INTAS
gants (No.97-804, No.00-221) and by DIETZ automation GmbH.}}
\author{Dietz J.O. \\
DIETZ automation GmbH,\\
Biedersbergweg 53\\
D66538 Neunkirchen / Saar\\
Germany\\
e-mail: DIETZautomation@aol.com \and Shashkov V.M. \\
Department of Numerical and Functional Analysis, \\
Faculty for Numerical Mathematics and Cybernetics, \\
Nizhny Novgorog University, 23/1 Gagarin avenue, \\
Nizhny Novgorod, 603600, Russia \\
e-mail: shashkov@unn.ac.ru \and Shashkov M.V. \\
Department of Differential Equations, \\
Institute for Applied Mathematics and Cybernetics, \\
10 Ulyanov street, Nizhny Novgorod, 603005, Russia\\
e-mail: shashkov@unn.ac.ru}
\date{}
\maketitle

\begin{abstract}
In this work we consider a class of quasilinear systems of differential
equations which allows to describe dynamics of electrohydraulic servovalves.
A method for fast identification of static and dynamic responses, by a
short-time experiment, is described.
\end{abstract}

\section{Introduction}

In 1950 Bill Moog invented a very reliable device to control flow in
hydraulic systems by electric signals. Since that time electrohydraulic
servovalves are key elements of various automatic control systems. Now it is
impossible to imagine rockets, aircrafts, satellites, certain machine tools
and control systems without servovalves. Such automatic systems require high
performance of the servovalve. That is why research and development of such
control systems as well as their maintenance require to test the dynamic and
static parameters of servovalves. Unfortunately up to now, test equipment
for servovalves has numerous of lacks. First of all, the test equipment
needs powerful hydraulic pump stations, which require water cooling, noise
insolation and frequent maintenance.\footnote{%
The power of pump stations reaches 100 kW and more.} Moreover, this
equipment requires large space and tubing for the distribution of the
hydraulic fluid. The installation of such test equipment takes time and
moving it is difficult. We must notice, that at the present time the
complete test equipment costs more than \$150,000 and a test of a servovalve
takes about 15 minutes. The cost may be increased many times if to take into
account that the different servovalves use different kinds of hydraulic
fluids (yellow oil, red oil, kerosene, skydrol, water and so on). Notice
that noise from the pump stations creates an additional problem for
hydraulic test equipment. The point is that the presence of noise
superimposes restrictions on the construction of buildings and the
laboratories where the test stand is located. Today, we know how to make a
quiet, compact, mobile and cheap hydraulic test stand.\footnote{%
The estimated expenses for the equipment are less \$25,000.}

In the present paper an approach for testing dynamic and static
characteristics of servovalves and principles of construction of modern test
equipment are described. The basic idea of the new approach is to divide the
test process into many test sub-processes, each of which takes a very short
interval of time. This idea implies the possibility to use the energy of a
hydraulic accumulator to produce a short subtest and, therefore, to make
obsolete powerful pump stations. Using this idea, we have developed a series
of mathematical methods, wrote computer programs and made computer
simulations for various servovalves. The experiments showed a possibility to
identify the static characteristic of a servovalve for 1 to 2 seconds of the
test time. To identify the dynamic characteristics we use a number of
sub-tests each of which requires just some milliseconds of the test time.
So, using a 4kW motor to fill a three-liter accumulator, the complete test
of a servovalve takes less than minute. It is necessary to note that the new
approach is based on the construction of a dynamic model of the servovalve
which is under test and this fact turns out to be extremely useful. The
knowledge of the dynamical system opens unlimited possibilities for
designers and developers of control systems. They can use Mathematica,
MathCad, MatLab, Simulink, LabView and other software for fast simulations
on a computer as well as for complex mathematical analysis.

\section{Statement of the problem}

A servovalve is a complex dynamical system and the complete description of
all its parameters, obviously, is not possible. However, most of engineers
and developers are interested in some basic dynamic and static
characteristics. To describe dynamics of servovalves, at the present time,
people usually consider amplitude-phase frequency responses. In this case,
obviously, it is assumed implicitly that dynamics of a servovalve can be
described well by a linear differential equation or by a linear transfer
function.\footnote{%
In a non-linear case the amplitude-phase response depends on amplitude and
form of the input signal and, therefore, does not make big sense.} In
practice, to describe the dynamic parameters, engineers usually use the
following system of differential equations of the second order (see for
instance \cite{Gu}, \cite{Ga}, \cite{Th1}).

\begin{equation}
A\frac{d^{2}x}{dt^{2}}+B\frac{dx}{dt}+Cx=u(t), \ \ \frac{dx}{dt}%
(0)=x(0)=0.  \label{lin2}
\end{equation}
Here, $u(t)$ is an input electrical signal, $x(t)$ is the output flow of a
liquid, $A$, $B$ and $C$ are some fixed parameters of the servovalve.

The static characteristic $x=F(u)$ of the servovalve is close to the linear
function $x=u/C$, however, in practice, there is the necessity to know the
static behaviour much better. In this case, by definition, to construct the
static response, it is necessary to know the output flow $x$ at any fixed
control $u$ \cite{Th2}. Here, of course, it is supposed, that at fixed
control the servovalve is a structurally stable dynamic system with a unique
stable equilibrium point. The such direct construction of the static
characteristic, obviously, requires time because it is necessary to check
the valve at all the points of the control range and each of the points
requires a delay for the transition time. In practice, to find the statics
approximately, people usually construct the Lissajous figure $(u(t),x(t))$,
where $u(t)$ is a very low-frequency harmonic control signal (usually less $%
0.01$Hz). But in any case it takes time.

In the present work, to take into account the nonlinearity of the static
characteristic, we shall consider the following quasilinear dynamic model of
a servovalve. 
\begin{equation}
A\frac{d^{2}x}{dt^{2}}+B\frac{dx}{dt}+f(x)=u(t), \ \ \frac{dx}{dt}%
(0)=x(0)=0.  \label{nonlin2}
\end{equation}
This system differs from the system (\ref{lin2}) by the presence of a
non-linear function $f(x)$. In general, such a nonlinearity can lead to the
appearance of a very complex chaotic behavior.\footnote{%
In 1918 Duffing found a chaos in a similar system with a cubic nonlinearity.
See, for instance, \cite{Guck}.} However, a real servovalve is a stable
dynamical system and its static response is close to a linear function.
Therefore we shall require that the function $f(x)$ is a monotone one and is
close to the function $Cx$. Moreover we shall assume, that the roots of the
characteristic equation $A\lambda ^{2}+B\lambda +C=0$ of the corresponded
linear system have negative real parts. These conditions guarantee a
stability of the system (\ref{nonlin2}) at any constant control $u(t)=const$
and, therefore, the static characteristic $x=F(u)$ of the servovalve is the
non-linear function $x=f^{-1}(u)$. The experiments made by the authors of
the present work have shown that the system (\ref{nonlin2}) describes well
both static and dynamic characteristics of servovalves. Therefore the main
problem to test a servovalve is reduced to the identification of the
parameters $A$, $B$ and the function $f(x)$ of the system (\ref{nonlin2}).
To solve this problem it is necessary to suggest a form of the test signal $%
u(t)$ and give an algorithm for a numerical analysis of $x(t)$ which
identifies $A$, $B$ and $f(x)$. Notice that we are interested in such
methods which give a good solution of the problem for a short enough time of
testing. Moreover, the methods should be structurally stable to the presence
of noise and measurement errors.

\section{Amplitude-phase frequency responses of a quasilinear system}

Let a solution $x(t)$ of the system (\ref{nonlin2}) correspond to the input
harmonic signal $u(t)=\sin (\omega t)$. After some transition time $\Delta t$
the output signal $x(t)$ will be a periodic function with the same frequency 
$\omega $. In this case, obviously, $x(t)$ can be represented by the
following Fourier series expansion: 
\[
x(t)=\sum_{k=0}^{\infty }R_{k}(\omega )\sin (k\omega t+\varphi _{k}(\omega
)). 
\]
For any $k$, the amplitude $R_{k}(\omega )$ and the initial phase $\varphi
_{k}(\omega )$ of the corresponding harmonic are calculated by the formulas: 
\begin{equation}
\begin{array}{l}
R_{k}(\omega )=\left| K_{k}(i\omega )\right| , \ \ \varphi
_{k}(\omega )=\arg \left( K_{k}(i\omega )\right) , \\ 
\\ 
K_{k}(i\omega )=\dfrac{\omega }{2\pi }\dint\limits_{\Delta t}^{\Delta t+2\pi
/\omega }x(t)e^{-ik\omega t}dt.
\end{array}
\label{K}
\end{equation}
Let us note that the amplitude $R_{1}(\omega )$ and the phase lag $\varphi
_{1}(\omega )$ of the first harmonic are known as the amplitude-phase
characteristics of the dynamic system. Namely these functions represent the
most of interest at the research and the development of automatical control
systems.

There are many different methods to construct the frequency responses $%
R_{1}(\omega ),$ $\varphi _{1}(\omega )$. The most simple method is to use
the formulas (\ref{K}) to calculate the functions directly. In this case, to
construct the frequency response, it is necessary to make measurements for a
larger series of frequencies $\omega =\omega _{1},\ldots ,\omega _{N},$
(usually $N>20$) and, then, to find the corresponding interpolation curves $%
R_{1}(\omega ),$ $\varphi _{1}(\omega )$.

The other method is to identify the parameters $A$, $B$ and function $f(x)$
of the quasilinear system (\ref{nonlin2}). In this case the test function $%
u(t)$ can be chosen from a wide class of functions (for example we can use
the step signal $u(t)=1(t)$). In order to identify $A$, $B$ and $f(x)$ the
following optimisation problem, for instance, can be solved: 
\begin{equation}
\min_{A,B,f}\dsum\limits_{k=1}^{N}\left( x(t_{k})-x_{A,B,f}(t_{k})\right)
^{2}.  \label{min1}
\end{equation}
Here $x(t_{k})$ are the results of measurement (samples) of the output
signal and $x_{A,B,f}(t_{k})$ is a solution of the quasilinear system (\ref
{nonlin2}) with parameters $A$, $B$, $f$ \ for the control function $u(t)$
at the moments $t_{k}\in (\tau _{1},\tau _{2})$. The minimization is
performed with respect to $A$, $B$ and $f$, where $f$\ is a function of a
finite-parametric family of functions $f_{\mu }$. So, the minimization is
performed with respect to parameters $A$, $B$ and $\mu $. The knowledge of
parameters $A$, $B$ and $f$, obviously, allows to construct the frequency
responses. This method allows to identify the parameters of the system for
an arbitrarily small time interval of measurement $(\tau _{1},\tau _{2})$.
However, the method requires a high precision in measurement of the output
signal $x(t_{k})$.

In our case to find the frequency responses $R_{1}(\omega ),$ $\varphi
_{1}(\omega )$ of servovalves it is expedient to use a method for
construction of dynamic characteristics which takes into account the
proximity of the system (\ref{nonlin2}) to the linear one (\ref{lin2}). In
this case, it is natural to suppose, that $R_{1}(\omega ),$ $\varphi
_{1}(\omega )$ are very similar to the frequency responses of the linear
system (\ref{lin2}) and our problem is reduced to search the optimum
parameters $A^{\ast }$, $B^{\ast }$ and $C^{\ast }$. Of course, $A^{\ast }$, 
$B^{\ast }$ and $C^{\ast }$ can be found as a solution of (\ref{min1}) but
we suggest an other way. It is well known, that the frequency responses of
the linear system (1) can be computed by the following formulas: 
\begin{equation}
\begin{array}{l}
R(\omega )=\left| K(i\omega )\right| , \ \ \varphi (\omega )=\arg
(K(i\omega )), \\ 
\\ 
K(s)=1\diagup \left( As^{2}+Bs+C\right) .
\end{array}
\label{Klin}
\end{equation}
Here $K(s)$ is the transfer function of the system (\ref{lin2}). In order to
find parameters $A^{\ast }$, $B^{\ast }$ and $C^{\ast }$ we suggest to solve
the following optimisation problem: 
\begin{equation}
\min_{A,B,C}\dsum\limits_{k=1}^{N}\left| K(i\omega _{k})-K_{1}(i\omega
_{k})\right| ^{2},  \label{min2}
\end{equation}
where $\omega _{1},\ldots ,\omega _{N}$ is a series of test frequencies and $%
K_{1}(i\omega _{1}),\ldots ,K_{1}(i\omega _{N})$ are values of the transfer
function which were computed by (\ref{K}). In contrast to the above
mentioned methods the last one does not require high precision in
measurement of the output signal $x(t_{k})$ and gives good results fast
enough because the method works for small series of experiments ($N\geq 2$).

\section{The DSS method}

In order to identify the static response we use the DSS method. The basic
idea of this method is the following. Let $x(t)$ be an output signal of the
dynamic system (\ref{nonlin2}) for a test signal $u(t)$. Then, obviously,
the following relation is fulfilled: 
\begin{equation}
f(x(t))=u(t)-A\frac{d^{2}x}{dt^{2}}(t)-B\frac{dx}{dt}(t).  \label{f}
\end{equation}
That means that the vector-function

\begin{equation}
(x(t),f(x(t)))=\left( x(t),u(t)-A\frac{d^{2}x}{dt^{2}}(t)-B\frac{dx}{dt}%
(t)\right) \in \Bbb{R}^{2}  \label{parst}
\end{equation}
is a parametric form of the static characteristic $u=f(x)$ (or $x=f^{-1}(u)$%
). We can consider that the term 
\begin{equation}
-A\dfrac{d^{2}x}{dt^{2}}(t)-B\dfrac{dx}{dt}(t)  \label{dcor}
\end{equation}
in (\ref{parst}) is a dynamical correction of the Lissajous figure $\left(
x(t),u(t)\right) $. Namely this term allows to construct the static
characteristic of the system for a short time, i.e. without a delay for the
relaxation of transition processes. However, to evaluate the dynamic
correction (\ref{dcor}), it is necessary to know parameters $A$ and $B$ of
the system (\ref{nonlin2}). Below we describe a method of successive
approximations to find $A$, $B$ and $f(x)$. First of all, we solve the
problem (\ref{min2}) to find parameters $A^{\ast }$, $B^{\ast }$ and $%
C^{\ast }$. The found parameters $A^{\ast }$, $B^{\ast }$ and the linear
function $C^{\ast }x$ is the first approximation $A_{0}$, $B_{0}$ and $%
f_{0}(x)$ of the corresponding parameters $A$, $B$ and function $f(x)$ of
the system (\ref{nonlin2}). To construct the next approximations $A_{n+1}$, $%
B_{n+1}$ and $f_{n+1}(x)$ ($n=0,1,2,\ldots $) we use the following inductive
rule. Let assume that the parameters $A_{n}$, $B_{n}$ and function $f_{n}$
are known and $x_{n}(t)$ is a solution of 
\begin{equation}
A_{n}\frac{d^{2}x_{n}}{dt^{2}}+B_{n}\frac{dx_{n}}{dt}+f_{n}(x)=u(t),\ \ 
\frac{dx}{dt}(0)=x(0)=0.  \label{nonlin2n}
\end{equation}
Then:

\begin{enumerate}
\item  Substituting the parameters $A_{n}$ and $B_{n}$ into the formula (\ref
{f}) instead of $A$ and $B$, we find the next approximation $f_{n+1}(x)$ of
the function $f(x)$.

\item  The next generation $A_{n+1}$ and $B_{n+1}$ of the parameters $A$ and 
$B$ can be found as a solution of the following optimisation: 
\begin{equation}
\min_{A_{n=1},B_{n+1}}\dsum\limits_{k=1}^{N}\left(
x(t_{k})-x_{A_{n+1},B_{n+1},f_{n+1}}(t_{k})\right) ^{2}.  \label{min3}
\end{equation}

Here $x(t_{k})$ are the results of measurement (samples) of output signal
and $x_{A_{n+1},B_{n+1},f_{n+1}}(t_{k})$ is a solution of the quasilinear
system (\ref{nonlin2}) with parameters $A_{n+1}$, $B_{n+1}$, $f_{n+1}$ \ for
the control function $u(t)$ at the moments $t_{k}\in (\tau _{1},\tau _{2})$.
Notice, the minimization is performed with respect to $A_{n+1}$ and $B_{n+1}$
only.
\end{enumerate}

The described steps should be repeated until the acceptable approximation
will be found. One more significance of our method is that, there are no big
needs to compute the derivatives $\dfrac{dx}{dt}(t_{k})$ and $\dfrac{d^{2}x}{%
dt^{2}}(t_{k})$ numerically or to use for that special electronics. The
point is that we can use the derivatives of the approximated solution $%
x_{n}(t)\equiv x_{A_{n},B_{n},f_{n}}(t)$ which was obtained on the $n^{th}$
step of the iteration process. In this case the $n^{th}$ approximation of
the static characteristic is given by the following parametric form: 
\begin{equation}
(x(t),f(x(t)))=\left( x(t),u(t)-A_{n}\frac{d^{2}x_{n}}{dt^{2}}(t)-B_{n}\frac{%
dx_{n}}{dt}(t)\right) .  \label{parstn}
\end{equation}
By (\ref{nonlin2n}) the formula (\ref{parstn}) takes the following simple
form:

\[
(x(t),f(x(t)))=(x(t),f_{n}(x_{n}(t))). 
\]
So, to construct the static characteristic we do not need to compute the
derivatives of the output signal $x(t)$. That fact is extremely important.
The point is that the measured signal $x(t)$ contains noise and the
measuring errors which lead to big problems at the differentiation.

In conclusion we must note that we did not investigate all questions about
convergence of our DSS process however the computer simulations demonstrated
fantastic results. Moreover, the acceptable results can be obtained even on
the first step of the DSS process if to use a low frequency input signal ($%
1-2$ Hz).

Notice also that the DSS method can be very simply generalised onto a
quasilinear systems of high order. For example, we tested our approach for
the following system of differential equations of the third order: 
\[
A\frac{d^{3}x}{dt^{3}}+B\frac{d^{2}x}{dt^{2}}+C\frac{dx}{dt}+f(x)=u(t),
\ \ \frac{d^{2}x}{dt^{2}}(0)=\frac{dx}{dt}(0)=x(0)=0. 
\]

\end{document}